\documentclass[11pt]{article}
\usepackage{amssymb}
\usepackage{amsmath}

\def\C{\mathbb C}                   
\def\K{\mathbb K}                   
\def\r{\mathbb R}                   
\def\fra{\mathfrak}                 
\newtheorem{defin}{Definition}

\newtheorem{propos}{Proposition}

\begin{document}

\pagestyle{myheadings}
\markboth
{{\sc Elements of the theory of induced representations \dots}}  
{{\sc O. Arratia and Mariano A. del Olmo}}        
\vspace{1cm}


\thispagestyle{empty}

\begin{center}

{\Large{\bf{Elements of the theory of induced representations}}}\footnote{Talk given by 
M.A.O. at the VIII Encuentro de Geometr\'{\i}a y F\'{\i}sica (Valencia, 
September 1998)}    

\medskip

{\Large{\bf{for quantum groups}}}     

\vskip1cm

{\sc O. Arratia $^1$ and Mariano A. del Olmo $^2$}
\vskip0.5cm

{ \it $^{1}$ Departamento de Matem\'atica Aplicada a la Ingenier\'{\i}a,  \\
Universidad de Valladolid. E-47011, Valladolid, Spain. 
\vskip0.3cm
$^{2}$ Departamento de  F\'{\i}sica Te\'orica,
 Universidad de Valladolid,  \\
 E-47011, Valladolid,  Spain}
\vskip0.15cm
 
\end{center}
\bigskip


\begin{abstract}
We analyze the elements  characterizing the theory of induced representations of
 Lie groups, in order to generalize it to quantum groups.
We emphasize the geometric and algebraic aspects of the theory, because they are
more suitable for generalization in the framework of Hopf algebras. As an
example, we present the induced representations of a quantum deformation of the
extended Galilei algebra in ($1+1$) dimensions.
\end{abstract}

\section{Introduction}

The theory of group representation plays a fundamental role in theoretical and
mathematical physics. Groups and physics are connected by the concept of symmetry,
i.e., any transformation on a system that keep
it invariant in some sense to make explicit. 
Induced representations are fundamental in theory of representations of Lie
groups and in many of their physical applications. If the  configuration space, $X$,
of a  physical system  is a  homogeneous space  of a certain Lie group of
transformations $G$, the (unitary) representations supported by this system can
be obtained  inducing from the  (unitary) representations of a Lie subgroup $K$ of
$G$, such that $X\simeq G/K$. So, there exists a deep relation between geometry and
representation theory of Lie groups via induced representations. 

One of the applications of quantum algebras and quantum groups is as generalized
symmetries of the symmetries associated to Lie groups. So, it looks natural to 
study the induced representations of quantum groups. Some attempts has been
made up to now. In Ref. \cite{PW}--\cite{ciccoli98} 
the  induction procedure of corepresentations  for quantum groups is
analyzed in a general way.  Corepresentations for some deformations of the 
Galilei group can be found in Ref.  \cite{m94,gkmt98}. Induced
representations of quantum groups has been studied in some particular cases:  the
one-dimensional quantum Galilei group in Ref. \cite{bgst98}, 
$sl_q(2)$ and $e_q(2)$ in Ref. \cite{ds94},  and  the quantum Heisenberg algebra in
Ref. \cite{es95}. 

 The aim of this paper is to analyze  and give an algebraic
formulation of the  elements of geometry as well as
of the theory of induced representations suitable to be generalized  to 
quantum algebras and groups. In this way we obtain an induction method of
representations for quantum groups \cite {arratia98}--\cite{arratia00}, whose 
 main advantage  is that it is not necessary the 
use of corepresentations or of Kirillov's theory \cite{Kir76}. 

 It is worthy to mention that the framework of modules and comodules is  the
appropriate one for our procedure, since it allows to profit  the algebraic
character of the structures involved. The concept of  pairing, which is the 
implementation of the idea of duality  to Hopf algebras (the algebraic structures
underlying  quantum groups and quantum algebras)  plays an important role,
since it shows  the equivalence of  modules/comodules and  spaces/co-spaces in the
theory of induced representations. By means of the duality starting from  regular
and  induced representations we get also coregular and coinduced ones. 

The main  conclusion of our analysis is that in the construction of induced 
representations is necessary to take in consideration two regular co-spaces:
one of them determines the equivariance condition which characterized the carrier
space of the induced representation, and a suitable  restriction of the action on
the other co-space gives the  induced representation itself.

Finally, to mention that in the example studied in section \ref{q-galileo}, the
quantum deformation of the extended Galilei $(1+1)$ algebra
$U_q(\overline{\fra{g}(1,1)})$, is fundamental to have  dual bases of it and its
dual $Fun_q(\overline{G(1,1)})$ to obtain the  coregular and coinduced
representations of $U_q(\overline{\fra{g}(1,1)})$.

\section{Algebraic preliminaries}

\subsection{Hopf algebras}
We present a brief review of the main ideas about 
Hopf algebras and that are relevant for our purpose (for more details see, for
instance, Ref. \cite{ChP}). 

\begin{defin}
An associative algebra with unit element over a commutative field,
$\mathbb{K}$ $(\r$ or $\C )$, is a triad $A=(V, m, \eta)$, where  $V$ is a linear
vector space and the maps
$ m:V\otimes V \longrightarrow  V$ and $\eta:  \mathbb{K} \longrightarrow V$ 
are linear and verify:
$$ \label{alg}
 m\circ(m \otimes \text{id})= m \circ (\text{id} \otimes m), \qquad
m \circ (\eta \otimes \text{id})= \text{id} = m 
\circ (\text{id}\otimes\eta).       
$$
\end{defin}
\begin{defin}
A coassociative coalgebra with  counit element over  a commutative field,
$\mathbb{K}$, is a triad $C=(V, \Delta, \epsilon)$, where $V$ is a linear
vector space and the maps coproduct, 
$
\Delta :V \longrightarrow  V\otimes V$,  and counit,
$\epsilon: V \longrightarrow \mathbb{K}$,
are  linear and verify
$$
 (\Delta \otimes \text{id}) \circ \Delta =(\text{id} \otimes \Delta)
\circ \Delta, \qquad
(\epsilon \otimes \text{id}) \circ  \Delta =\text{id} =
   (\text{id} \otimes \epsilon) \circ \Delta.        
$$
\end{defin}    
In general,  coalgebras have a behaviour ``dual'' of that of  algebras. 
The coproduct of an element $c\in V$ can be written as
$\Delta(c)= \sum_{(c)} c_{(1)} \otimes c_{(2)}$.
 \begin{defin}  
A bialgebra is a pair $B\equiv(A,C)$ made up of an algebra $A$
 and a coalgebra $C$, both with the same underlying vector space 
$V$, such that  $\Delta$  and $\epsilon$  are morphisms of $A$.
\end{defin}
\begin{defin}
A bialgebra $H=(V,m,\eta,\Delta,\epsilon)$
 is said to be a Hopf algebra if there exists an antipode, i.e., a bijective linear
map  $S:V \rightarrow V$ verifying
$$
m \circ (S \otimes \text{id}) \circ \Delta= \eta \circ \epsilon=
m \circ ( \text{id} \otimes S) \circ \Delta.
$$
\end{defin}   
 It can be proved that the antipode is an anti-automorphism of algebras and an
anti-coautomorphism of coalgebras, and if it exists then  is unique.  

As examples of Hopf algebras we can mention: finite group
algebras, the algebra of functions on a finite Lie group,
and enveloping algebras of Lie algebras. All of them have the
property of commutativity or  cocommutativity, i.e., 
$m \circ \tau= m$  or $\tau \circ \Delta =\Delta$, with $\tau$ the twist operator,
i.e., $\tau(a\otimes b)=b\otimes a,\ a,b \in V$.
All of them are non-deformed  or ``classical".

\subsection{Quantum algebras and quantum groups}

Quantum groups and quantum algebras are  Hopf algebras which
are  neither commutative nor cocommutative. The usual
definition of quantum algebra (in the sense of Drinfel'd \cite{drinfeld86}
and Jimbo \cite{jimbo85}) is as follows.
\begin{defin}
Let ${\cal U}({\fra g})$  be the universal enveloping algebra of a Lie algebra
${\fra g}$. It is a ``classical" Hopf algebra  with coproduct, counit and antipode
defined by
$$
\begin{array}{ccccc}
&\Delta (X)=1\otimes X+ X\otimes 1, \qquad 
&\ \ \Delta (1)= 1\otimes 1;\\[0.2cm]
&\epsilon (X)= 0,\quad \epsilon (1)= 1; \quad
&S (X)= -X ; 
\end{array}\qquad  X \in {\fra g} .
$$
The extension to the remaining elements of ${\cal U}({\fra g})$ is made by
linearity. 
\end{defin}

A quantization or deformation of ${\cal U}({\fra g})$ is obtained by
means of a deformed Hopf structure on the associative algebra of formal power
series in $z$ and coefficients in ${\cal U}({\fra g})$, 
${\cal U}_z({\fra g})\equiv {\cal U}({\fra g})\hat\otimes \C[[z]]$,  such
that ${\cal U}_z({\fra g})/z{\cal U}_z({\fra g})\simeq {\cal U}({\fra g})$,
i.e., ${\cal U}_z({\fra g})\to  {\cal U}({\fra g})$ when $z\to 0$ 
(the equivalence is in the sense of Hopf algebras). 

On the other hand, there are several approaches for
quantum groups \cite{fadeev,woronowicz}. Let $G$ be a finite
dimensional Lie group and  ${\fra g}$ its Lie algebra. Let us consider the
commutative and associative algebra of smooth functions of $G$ on $\C$, $Fun(G)$, 
with the usual product of functions (i.e., $(fg)(x)=f(x)g(x), \ f,g
\in  Fun (G),\ x,y \in G$). It has a Hopf structure  (which is cocommutative if
and only if $G$ is abelian) given by
$$
(\Delta (f)) (x,y)= f (xy),\qquad  \epsilon (f)= f(e), \qquad
(S(f))(x)= f(x^{-1}) ,
$$
where $e$ is the identity of $G$.   

Note that, in general, $Fun (G)\otimes Fun (G) \subseteq  Fun (G\times G)$.
When the group is finite the equality is strict, but if $G$ is not
a finite group  $\Delta (f)$ may not belong to 
$Fun (G)\otimes Fun(G)$. This problem can be solved by an
adequate restriction of the space $Fun (G)$.

After deformation the above commutative Hopf algebra becomes
non-commutative $Fun_q (G)$. 
In section \ref{q-galileo} we  display the quantum Galilei algebra
$U_q(\overline{\fra{g}(1,1)})$ as well as of the quantum Galilei group
$F_q(\overline{G(1,1)})$.

\subsection{Pairing of Hopf algebras}

In the category of linear
vector spaces ($V$) over $\K$ the dual object of $V$ is defined
as the vector space of its linear forms, i.e., 
$V^*= {\cal L}(V, \K)$.

If $(V,m,\eta)$ is a finite algebra its  dual object  $(V^*, m^*,\eta^*)$ is
define to be  a coalgebra. Reciprocally, the dual of a coalgebra is an algebra.
However, the category of bialgebras and Hopf algebras is self-dual.

When $V$ is infinite dimensional,  $(V \otimes V)^*$ and
$V^* \otimes V^*$ are not isomorphic (in fact, there is only an
inclusion of $V^* \otimes V^*$ into $(V \otimes V)^*$).
Hence, the dual of a product map, 
$m: V \otimes V \rightarrow V$, 
does not give rise to a coproduct in a
natural way, but the concept of pairing \cite{ChP} solves this difficult.

\begin{defin}
A pairing between two Hopf algebras, $H$ and $H'$, is a bilinear map
$\langle \, \cdot \, ,  \, \cdot \, \rangle:
H \times H' \rightarrow \K$ that verifies the following
properties:
$$\label{pairing}
\begin{array}{ll}
\langle h , m'(h'\otimes k')\rangle= \langle \Delta(h) , h' \otimes k' \rangle, 
&  \quad\langle h, 1_{H'} \rangle= \epsilon(h),\qquad
 \epsilon'(h')= \langle  1_{H}, h' \rangle ,\\[2mm]
\ \langle h\otimes k , \Delta'(h')\rangle=
\langle m(h\otimes k) , h' \rangle, & 
\langle h, S'(h') \rangle =     \langle S(h), h' \rangle.  
\end{array}
$$
\end{defin}
The pairing is left (right) non-degenerate if
$[\langle h, h' \rangle = 0,  \forall h' \in H'] \Rightarrow h=0$
($[\langle h, h' \rangle = 0, \ \forall h \in H] \Rightarrow h'=0$).  We
say that a pairing is non-degenerate if it is simultaneously left and
right  non-degenerate. A non-degenerate pairing
allows us to restrict (in the infinite dimensional case) the dual $H^*$
to a more ``manageable'' subspace isomorphic to $H'$ via
$\widetilde{h'}(h)= \langle h , h' \rangle$.

 Incidentally, $Fun (G)$ ($Fun_q (G)$) is the Hopf algebra dual of  
${\cal U}({\fra g})$ (${\cal U}_q({\fra g})$).

\subsection{Star structures on  Hopf algebras}\label{starhopf}

Complex numbers solve some problems that appear
using  real numbers, but the results of the physical measures are real
numbers. So, it is pertinent to introduce ``real forms''(star structures) on
the complex spaces  we are using. 
\begin{defin}\label{stardefin}   
A  star structure on a  Hopf algebra,  $H$, is a map
$*:H\rightarrow H$  such that 
 \begin{equation} \label{axe}
 \begin{array}{ccc}
 h^{**}= h, & (\lambda h + \mu g)^*= \bar{\lambda} h^* + \bar{\mu} g^*,
& (hg)^*= g^*h^*,\\[0.2cm]
\quad \Delta(h^*)=\Delta(h)^*, & \epsilon(h^*)= \overline{\epsilon(h)},&  1^*=1,
\end{array}
\end{equation} 
where $h,g \in H$, $\lambda, \mu \in \C$ and  the  bar denotes complex conjugate.
\end{defin}
The axioms displayed in expression (\ref{axe}) show that $*$ is an antimorphism of
algebras and a  morphism of coalgebras, which is the usual choice 
\cite{ChP,Maj95a}.  
 
In the set of  $\mathbb{C}$--valued functions on the space  $X$
the star structure is defined by $f^*(x)= \overline{f(x)}$.
If $G$ is a finite group, the  star structure to be considered on the 
algebra $\mathbb{C}[G]$ is given by $g^*= g^{-1}$, $\forall g \in G$.
When  $G$ is a Lie group, the canonical star structure on the
complexification of the  enveloping algebra $U(\mathfrak{g})\otimes \mathbb{C}$ 
is defined for the generators $X$ of $\mathfrak{g}$ as $X^*= -X$, and it is 
extended to all the remaining elements taking into account that $*$ is a semilinear
antimorphism of algebras. Note that  with Def. \ref{stardefin}, 
$S\circ *= S^{-1}\circ *$. 
 
Given a pair of Hopf algebras equipped with a non-degenerate pairing
$(H, H', \langle \cdot , \cdot \rangle)$, the star structure on $H$ can be
translated to $H'$ by means of $\langle \varphi^*, h \rangle =
\overline{\langle \varphi, S(h)^* \rangle}$.

\subsection{Modules and comodules}\label{moduloscomodulos}

The concept of module appears  from that of linear vector space  substituting the 
field of scalars by a ring, and  dualizing modules  we get comodules.
 
\begin{defin}  Let $(V, \alpha, A)$ a triad made of an associative
$\mathbb{K}$--algebra, $A$, with unit element, a linear vector $\mathbb{K}$--space,
$V$, and a linear map $\alpha: A \otimes_{\mathbb{K}} V \rightarrow V$, that we
call action, and denote by $a \triangleright v = \alpha(a\otimes v)$ over the 
decomposable elements of the tensor product. We will say  that 
$(V,\alpha, A)$ (or $(V, \triangleright, A)$) is a left $A$--module if  verifies:
$a \triangleright (b\triangleright v)= (ab) \triangleright v, 
\ 1 \triangleright v= v, \ 
 \forall a, b \in A,  \forall v \in V$.
\end{defin}
We only consider  ``left  objects''  but, obviously, 
right versions can be defined  in a similar way. In the regular
$A$--modules  $(A, \triangleright, A)$ and  $(A, \triangleleft, A)$ the actions are
obtained via the product of the  algebra by $a'\triangleright a=a'a$ and
$a\triangleleft a'=aa'$, respectively. The morphisms between two $A$--modules $(V,
\alpha, A)$ and 
$(V',\alpha, A)$ are the linear maps $f: V\to V'$ such that are equivariant respect
to the action, i.e., $f(a'\triangleright v)=a\triangleright ' f(v),\ 
\forall a\in A,\forall v\in V$.
 \begin{defin}
 Let $(V, \beta, C)$ be a triad made of a coassociative  $\mathbb{K}$--coalgebra
with counit, $C$, a linear vector $\mathbb{K}$--space,
$V$, and a linear map  $\beta: V \rightarrow C \otimes_{\mathbb{K}} V$,  that we
call coaction, and denote by $v \blacktriangleleft= \beta(v)  
= v^{(1)} \otimes v^{(2)}$.
We will say  that $(V, \beta, C)$ (or $(V, \blacktriangleleft, C)$) is a left
$C$--comodule if 
$$ \label{comodulo}
 {v^{(1)}}_{(1)} \otimes {v^{(1)}}_{(2)} \otimes {v^{(2)}}=
 {v^{(1)}} \otimes v^{(2)(1)} \otimes v^{(2)(2)}, \qquad
 \epsilon(v^{(1)}) v^{(2)}= v, \quad \forall v \in V.
$$
\end{defin}
The  regular right comodule  $(C, \blacktriangleleft, C)$  (left comodule
$(C, \blacktriangleright, C)$) has defined its coaction in terms of the
coproduct on $C$:
$ c \blacktriangleleft = c_{(1)} \otimes c_{(2)}$   
($\blacktriangleright c = c_{(1)} \otimes c_{(2)}$).
The morphisms between two $C$--comodules, $(V, \blacktriangleleft ,C)$ and 
$(V', \blacktriangleleft ',C)$, are linear maps 
 $f:V\rightarrow V'$ verifying
$v^{(1)} \otimes f(v^{(2)})=f(v)^{(1)'}\otimes  f(v)^{(2)'},
 \  \forall v \in V$.

The definition of  dual modules and  dual  comodules is possible when they are
building up over a Hopf algebra $H$.
Note that if $(V, \triangleright, H)$ is a left $H$--module then
the dual space $V^*$ is equipped in a natural way with a structure of right
$H$--module by
$\langle f \triangleleft h, v  \rangle =
\langle f ,h \triangleright  v  \rangle , \ f\in V^*, v \in V,   h \in H$.
When $(V, \blacktriangleleft, H)$ is a left $H$--comodule then $V^*$ is a right 
$H$--comodule by means of
$\langle \blacktriangleright f,  v\otimes \varphi  \rangle = 
\langle \varphi \otimes f ,  v \blacktriangleleft   \rangle , \ f\in V^*,v \in V, 
\varphi \in H^*$.  
The module dual of $(V, \triangleright, H)$ is obtained making the 
pull--back with respect  to the antipode $S$ of $H$, i.e., 
$h\triangleleft v =S(h)\triangleright v$. 
The comodule dual of $(V, \blacktriangleleft, H)$ is the push--out with respect to
$S$, i.e.,  $v\blacktriangleleft  =S(v^{(1)})\otimes v^{(2)}\triangleright v$. 

A representation of the algebra $A$ on a linear vector space $V$ is an algebra
morphism
 $\rho: A \rightarrow \text{End}(V)$.
The expression $ \rho(a)(v)= a \triangleright v$
 establishes a one-to-one correspondence between  representations of 
 $A$ and left $A$--modules, in such a way that any concept defined for
representations has its analogue in the language of modules. 

When a bialgebra acts or coacts in a linear vector space equipped with another
structure like algebra, coalgebra or bialgebra some compatibility relations are
necessary 
\cite{Maj95a,Mol}.
In the following  definitions $B$ and $B'$ denote bialgebras,
$H$ and $H'$ Hopf algebras, $A$ an algebra and $C$ a coalgebra.
\begin{defin}
The left $B$--module $(A, \triangleright, B)$ is an algebra-module
if $m_A$ and $\eta_A$ are morphisms of $B$--modules, i.e., if  
$ b\triangleright( aa')=(b_{(1)}\triangleright a)(b_{(2)}\triangleright a'),
$ and $ b \triangleright 1 = \epsilon(b) 1, \ \forall b \in B, \, \forall a, a' \in
A$.
\end{defin}
\begin{defin}  The left $B$--module  $(C, \triangleright, B)$ is a 
coalgebra--module if
$\Delta_C$ and $\epsilon_C$ are $B$--module morphisms, i.e., if 
for any  $b,c \in B$ it is accomplished that
$
(b\triangleright c)_{(1)} \otimes (b \triangleright c)_{(2)} =
    (b_{(1)}\triangleright c_{(1)}) \otimes (b_{(2)}\triangleright
c_{(2)})$ and $\epsilon_{C}(b \triangleright c) = \epsilon_{B}(b) \epsilon_{C}(c)$.
\end{defin}
By dualizing  the actions two new structures over a left $B$--comodule 
$(C,\blacktriangleleft, B)$ are obtained: coalgebra--comodule and 
algebra--comodule.
When the object where one acts is a bialgebra  other two structures are
obtained.
\begin{defin}
It is said that $(B', \triangleright, B)$  [$(B', \blacktriangleleft, B)$] is a left
bialgebra--[co]mod\-ule if is simultaneously an  algebra--[co]module and a
coalgebra--[co]module.
\end{defin}

\section{Spaces and co-spaces}\label{spacesyco-spaces}

The formalization of the idea of geometry is based on two objects: an space,  
$X$, and a group of transformations,  $G$, acting on it.
We understand for (right) action  of  $G$ on $X$
an external composition law $\alpha: X \times G \rightarrow X ,\
(\alpha(x, g)\equiv x \triangleleft g)$, verifying
$ (x \triangleleft g) \triangleleft g'= x \triangleleft (gg')$ and 
$ x \triangleleft 1= x$. 

So, we consider a geometry as a triad $(X, \triangleleft, G)$ constituted by an
space, $X$, a group, $G$, and a (right) action, $\triangleleft$,
 of $G$ on $X$. It is said that $X$ is a $G$--space or also a
$G$--module. 
For instance, the group $G$ with 
the regular action, $g' \triangleleft g=g'g, \ g,g \in G$, is a $G$--module. 
The morphisms between the  objects
$(X, \triangleleft, G)$ and $(X', \triangleleft', G)$ are  $G$--equivariant
maps $f:X\rightarrow X'$, i.e., $f(x \triangleleft g) =f(x) \triangleleft' g, 
\ \ \forall x \in X, \;  \forall g \in G$. 

The action of $G$ can be  extended in a natural way to other
objects defined in terms of $X$. For instance, the set $F(X,X')$ of 
maps of $X$ on $X'$ can be equipped with a (left)
action  defined by $(g \triangleright f)(x) =f(x \triangleleft g)$.
 
The fact to study the geometry in terms of other triads associated to the original
one is interesting for two reasons: some questions about 
$(X, \triangleleft, G)$ can be solved easier in the new object, and
 the characteristics of the new triad can allow a natural
generalization of the original geometric concepts.

In this way, we  replace the manifold $X$  by the algebra of ${\cal C}^\infty$
$\C$--valued functions on $X$, as well as the Lie group $G$ by the enveloping 
algebra $U(\mathfrak{g})$ of its Lie algebra $\mathfrak{g}$. Since  $({ F}(X),
\triangleright, U(\mathfrak{g}))$ is an algebra--module over the Hopf algebra
$U(\mathfrak{g})$, we can generalize the concept of $G$--space in algebraic terms.
\begin{defin}\label{h co-spaces}
Let $H$ be a Hopf algebra. A left [right] $H$--co-space is an algebra module 
$(A, \triangleright, H)$ [$(A, \triangleleft, H)$].
\end{defin}

The morphisms among \mbox{$H$--co-spaces} are the morphisms of
$H$--modules and the concepts of subco-space or quotient co-space are
equivalent to  subalge\-bra--module  or quotient algebra--module, respectively.
The term {\em co-space} is not usual, but  we have
adopted this word instead of {\em space} to stress the dual character of $A$ as way
of describing the initial ``geometric object''.
The crucial point in the above definition is that no condition about the
commutativity of the product in $A$ is required. Hence, it allows to write
in a geometric language properties of algebras which are not the ring of coordinates
 over a space. 

As examples we mention the following ones. Let $(X, \triangleleft,G)$ be a
$G$--space, with $G$ a finite group. The action of $G$  on the
$\mathbb{K}$--valued functions on $X$ by $(g \triangleright f)(x)= f(x
\triangleleft g),$ determines the $\mathbb{K}[G]$--module $( F(X), \triangleright,
\mathbb{K}[G])$, where $\mathbb{K}[G]$ is the group algebra. The right $G$--space
$(M,
\triangleleft, G)$, with $M$  a smooth manifold and $G$ a Lie group, has
associated the
$U(\mathfrak{g})$--co-space $(F(M), \triangleright, U(\mathfrak{g}) )$, where the
action is defined by 
$1 \triangleright f=f, \ (X \triangleright f)(p)= \frac{d}{dt} \big\vert_{t=0} 
f(p \triangleleft e^{tX})$, with $X \in \mathfrak{g}$.
If $H$ is a finite Hopf algebra by  dualizing the regular action we get the module
$(H^*, \triangleright, H)$ with the action given by $h \triangleright \varphi =
\langle h , \varphi_{(2)} \rangle \varphi_{(1)},  \ h \in H, \ \varphi \in H^*$.
In this example the  module is also a $H$--co-space.

In the non-finite case it is necessary to consider a pair of algebras with a
non-degenerate pairing $(H,  H', \langle \cdot, \cdot \rangle)$. It allows, via 
dualization of the regular actions, to obtain the regular $H$--co-spaces
\mbox{ $(H', \triangleright, H)$} and \mbox{ $(H', \triangleleft, H)$}. 
 
In the context of algebra of functions a subspace will be a quotient co-space and a
quotient space  a subco-space. 

In the  co-spaces $({ F}(X), \triangleright, \mathbb{C}[G])$ and
$({ F}(M), \triangleright, U(\mathfrak{g}))$,
the invariance of a function $f$ under the action means that
$ g \triangleright f = f, \ \forall g \in G,$ and 
$ X \triangleright f = 0, \ \forall X \in \mathfrak{g}$, respectively.
Both expressions are summarized in the general case $(A, \triangleright, H)$
saying that $a\in A$ is invariant under the action of $H$ if
$h \triangleright a =\epsilon(h) a, \ \forall h  \in H$.

As example of algebraic description  we will analyze 
the compatibility between the action on the  co-space $(A, \triangleright, H)$
and the star structures defined on $A$ and $H$. From 
 $({ F}(M), \triangleright, U(\mathfrak{g})\otimes \mathbb{C})$, now $F(M)$ 
denotes the space of ${\cal C}^\infty$ $\mathbb{C}$--valued functions on $M$,
equipped with star structures defined in section \ref{starhopf},  one gets:
$ (X \triangleright f)^*(p)=
\overline{\frac{d}{dt} \big\vert_{t=0} f(p \triangleleft e^{tX})}
=\frac{d}{dt} \big\vert_{t=0} f^*(p \triangleleft e^{tX})= 
(X \triangleright f^*)(p)$,
with $X \in \mathfrak{g}$, $f \in { F}(M)$ and $p \in  M$.
When $Z= X + i Y \in \mathfrak{g\otimes \mathbb{C}}$ we have
 $(Z \triangleright f)^* =(X \triangleright f)^*- 
i (Y \triangleright f)^* =  X \triangleright f^*- i  Y \triangleright f^*  
=- Z^* \triangleright f^*$.
In the  co-space $({F}(X), \triangleright, \mathbb{C}[G])$
with the star structures displayed in section \ref{starhopf}
we obtain the identity 
$(g \triangleright f)^*= g \triangleright f^*$.
The last two expressions  can be rewritten using the
antipode as $(Z \triangleright f)^* = S(Z)^* \triangleright f^*$ and 
$(g \triangleright f)^* =S(g)^* \triangleright f^*$, respectively,
which give  a compatibility relation for the general case:
$(h\triangleright a)^* =S(h)^* \triangleright a^*$.

Other  important concept to be generalized is the  measure on a space. Its 
algebraic expression is immediate, just like  to note that in a
space $X$ with a measure $\mu$  a functional on the set of
measurable functions is  defined by
$$
 I_\mu(f)= \int_X f(x) \, d\mu(x).
$$ 
In fact Riesz's theorem guarantees, under hypotheses very few
restrictive, the existence of a one-to-one correspondence between measures on 
 $X$ and positive linear functionals.

If the $H$--module  $(A, \triangleleft, H)$ is finite we can 
consider the dual module, $(A^*, \triangleright, H)$, where
 the integral $I \in A^*$
is left invariant if 
$h \triangleright I= \epsilon(h) I$.
\begin{defin} 
A left invariant integral on the  co-space  $(A,\triangleleft, H)$
is a linear form, $I:A \rightarrow \mathbb{K}$, that verifies 
$\langle I, a \triangleleft h \rangle 
       =  \epsilon(h) \langle I, a \rangle ,\  
\forall h \in H,\   \forall a \in A$.
The integral $I$ is normalized if $\langle I, 1_A \rangle =1$.
\end{defin}
For the  regular co-space $(H', \triangleleft, H)$, presented below,  the
condition that determines the invariance of the integral
$I:H' \rightarrow \mathbb{K}$ can be reformulated using the coproduct of
$H'$ as  $\langle h \otimes I, \Delta h' \rangle
=\epsilon(h) \langle I, h' \rangle,
\ \forall h \in H,  \ \forall h'\in H'$.
In this  case it is possible to demonstrate that such integral is unique up to
a constant. 

The definition of
homogeneous space is not easy to generalize and in the
literature there are some generalizations non-equivalent 
\cite{Bon96a,Pod87}.

It is worthy to mention that the use of the duality supplies  several alternatives
for the same object. So,  associated to the  $H$--co-space $(A, \triangleright,
H)$ there is another module,  $(A^*, \triangleleft, H)$, and two comodules,
$(A, \blacktriangleright, H^*)$ and $(A^*, \blacktriangleleft, H^*)$.


\section{Induced representations}\label{representacionesinducidas}

Let us start rewriting  the method of induced representations from an algebraic 
point of view in terms  of  modules \cite{Dix} and giving its dual version.
\begin{defin} Let $A'$ be a subalgebra of an algebra $A$, which
can be considered simultaneously as a right $A'$--module  and a left $A$--module 
for the regular action, and 
$(V, \triangleright, A')$ a left $A'$--module. Then, $A \otimes_{A'} V$
equipped  with a left  action of $A$ is  a left  \mbox{$A$--module}
$(V^\downarrow, \triangleright, A)$ called 
\mbox{$A$--module} induced by $(V, \triangleright, A')$. If $\rho$ and  $\pi$ are
the representations associated to $(V, \triangleright, A')$ and $(V^\downarrow,
\triangleright, A)$, respectively,  it is said that $\pi$ is the
representation of $A$ induced by $\rho$.
\end{defin}
\begin{defin}
Let $A'$ be a subalgebra of an algebra $A$, which can be considered 
simultaneously as a left $A'$--module and a right $A$--module for the
regular action, and $(V, \triangleright, A')$ 
a left $A'$--module.  Then the space $V^\uparrow =\text{Hom}_{A'}(A, V)$ 
equipped with  the action given by
$(a \triangleright f)(b)= f(ba), \  \forall a,b \in A, \;  \forall f \in
\text{Hom}_{A'}(A, V),$ is a left $A$--module
$(V^\uparrow, \triangleright, A)$,  called $A$--module  coinduced by 
$(V,\triangleright, A')$. If $\rho$ and $\pi$ are the representations
associated to $(V,\triangleright, A')$  and $(V^\uparrow, \triangleright,
A)$, respectively, it is said that  $\pi$ is the representation of $A$ coinduced
by $\rho$.
\end{defin}
The above procedures of  induction and coinduction appear in a natural way in this
context since they are extensions of scalars in a module \cite{Hig}.
For infinite dimensional algebras the coinduced representations 
have non-countable dimension, this difficult is
avoided  using non-degenerate
dual forms.

 Let $\langle \, \cdot \, , \, \cdot \, \rangle$  be a non-degenerate pairing
between the Hopf algebras $H$ and $H'$. If  $K$ is a subalgebra of $H$ and
$(V,\triangleright, K)$ a $K$--module, the carrier space
$V^\uparrow$ of the induced  module is the subspace of
$H'\otimes V$ with elements $f$  verifying
\begin{equation}\label{pareado}
\langle f, kh \rangle=  k \triangleright \langle f, h \rangle,
\qquad \forall k \in K,\ \forall h \in H.
 \end{equation} 
The pairing of  expression (\ref{pareado}) is $V$--valued and given by
$ \langle \varphi \otimes v, h \rangle =\langle \varphi , h \rangle v$, 
with $h \in H,  \varphi \in H',  v \in V$.
Finally, the action $h \triangleright f$ in the  coinduced module
is defined by
$$
\langle h \triangleright f, h' \rangle = \langle f, h'h \rangle,
 \qquad  \forall  h' \in H.  
$$

Let us suppose a  real Lie group $G$ and its Lie algebra $\mathfrak{g}$.  For any
unitary representation of $G$ on a complex Hilbert space $\cal H$ there are a
representation of  $U(\mathfrak{g}\otimes \mathbb{C})$ on the space of ${\cal
C}^\infty$ vector fields  of $\cal H$, and other of 
$U(\mathfrak{g}\otimes \mathbb{C})$ on the  space ${\cal H}^\infty$ of
distribution vectors of $\cal H$ \cite{Dix}.

\subsection{Group representations}
\label{representacionesdegrupos}

In this section we show the importance of the associative algebras in
the construction of group representations and, in particular, 
induced representations.

In the context of section \ref{spacesyco-spaces}, 
it is well known that any finite group 
$G$ has associated in a natural way the group algebra $\mathbb{K}[G]$, whose 
elements can be interpreted either as formal  linear combinations of elements of
$G$ with coefficients on $\mathbb{K}$ or as $\mathbb{K}$-valued functions on $G$. 
Moreover, the category of representations on vector space  over
$\mathbb{K}$ is equivalent to that of  the  modules over $\mathbb{K}[G]$
 \cite{Kir76}. 

If  $K$ is a subgroup of a finite group $G$ and  $(V, \rho)$ a
representation of $K$ one can consider the algebras $\mathbb{K}[K]$ and
$\mathbb{K}[G]$ together with the  module  $(V, \triangleright, \mathbb{K}[K])$
defined by $k\triangleright v = \rho(k)(v), \ k \in K, \; v \in V$.
The  carrier space  of the  induced representation,
$\text{Hom}_{\mathbb{K}[K]}(\mathbb{K}[G], V) $, is the set of linear maps $F$
defined in $\mathbb{K}[G]$  verifying $F(k g)=  k\triangleright F(g), \ \forall k
\in \mathbb{K}[K], \;
\forall g \in \mathbb{K}[G]$.
Since  $G$ is a basis of $\mathbb{K}[G]$  only it is necessary to know the values
of  $F$ on the elements of $G$. This allows to characterize the carrier space,
$F_\rho(G, V)$, of  the  induced representation as the subspace of $V$--valued
functions on $G$  verifying the  ``equivariance condition''
$$ 
f(kg) =\rho(k)(f(g), \qquad f= F|_{G},\; \forall k \in K, \;  \forall g \in G. 
$$
Note that the equivariance
condition is completely  determined by $K$, $G$ and $(V,\rho)$. The 
induced representation  $\rho^\uparrow$ is given by
$$ 
[\rho^\uparrow(g) f](g')= (g\triangleright f)(g')= f(g'g). 
$$

When the  group $G$ is not finite,  the  equivalence between representations and
$\mathbb{K}[G]$--modules is maintained, but the algebra is too complicated  and
usually other algebras are introduced \cite{Kir76}.

Let us suppose now that $G$ and $K$ are connected Lie groups and, for the sake of
simplicity,  $V$ is finite dimensional. Then, the algebras to be consider are
$U(\mathfrak{g})$ and $U(\mathfrak{k})$ together with the  module  
$(V, \triangleright, U(\mathfrak{k}))$ defined by
$X \triangleright v =\frac{d}{dt}\big\vert_{t=0} \rho(e^{tX})(v),
\ \forall X \in \mathfrak{k}$. 
The  carrier space of the  induced representation is the set of the linear 
$V$-valued functions on $U(\mathfrak{g})$ verifying
$$
  F(X Y) = X \triangleright F(Y), \qquad  \forall X \in U(\mathfrak{k}), \;
Y \in U(\mathfrak{g}). 
$$

\subsection{Induced corepresentations}

Taking into account the concept of duality  the knowledge of  the category
of representations of the Hopf algebra $H$  is equivalent to that of the
category  of corepresentations of the  dual of $H$ (at least if
$\text{dim}\; H < \infty$). So, it looks natural to ask ourselves for the
analogue of the induction algorithm in the category of comodules.
  In the following we discuss the  finite dimensional case. 

Let us consider the subalgebra $A'$ of the algebra
$A$ with canonical injection 
$  i: A' \longrightarrow A$, and the  module  $(V, \triangleright, A')$.
As we seen below  the induction method   originates
the $A$--module $ (V^\uparrow , \triangleright, A)$ with action 
$(a \triangleright F)(b)= F(ba)
$ and carrier  space $V^\uparrow =   \text{Hom}_{A'}(A, V) 
\subset\text{Hom}_{\mathbb{K}}(A,V)$, whose elements $F$  verify 
$F(a'a) = a' \triangleright F(a)$,  $\forall (a',a) \in A' \times A$.
The dualization gives the coalgebras
$C' ={A'}^{*}$ and $ C =A^*$
with canonical subjection 
 $\pi =i^*: C \longrightarrow C'$ and the comodules
$(V, \blacktriangleright, C'), \ (V^\uparrow, \blacktriangleleft, C).
$  The last comodule  is said to be the comodule  induced by  the former one.
Note that the  carrier  space is  the same in every corresponding pair
(module$/$comodule).

The above procedure can be  reformulated in terms of comodules. For that
it is necessary to rewrite the condition that characterizes  
the carrier space using coactions and  determine the coaction  on the induced
comodule.

The elements of $ V^\uparrow$ are of the form
$F ={\sum_{i \in I} } v_i \otimes \varphi_i$,
with  $I$ a set of indices,  $v_i \in V$ and $\varphi_i\in C$.
For all $f \in V^*$, $a'\in A'$ and $a \in A$ we have that
$ \langle F(a'a), f \rangle =\langle (\text{id} \otimes \Delta) F, 
f \otimes  a' \otimes a \rangle$.
Using the equivariance condition we get, for all $f \in V^*$, $a'\in A',\ a \in
A$, that
$\langle F(a'a), f \rangle   =\langle  (\blacktriangleright \otimes 
\text{id}) F, f \otimes a'\otimes a\rangle$.
Comparing both expressions  the equivariance condition is equivalent to 
$$ 
(\text{id} \otimes L ) F =(\beta \otimes \text{id})F,
$$
where $L =(\pi \otimes \text{id}) \circ \Delta$ and $\beta$ is the coaction on
the $C'$--comodule.

Noting that 
$\blacktriangleright F \in (V \otimes C)\otimes C$    
we can obtain  the coaction on $F$. So, 
$\langle (f \otimes a) \otimes b, \blacktriangleright F \rangle =
\langle f \otimes a \otimes b, (\text{id } \otimes \Delta ) F \rangle$.
From this result we deduce that the induced coaction is given by 
$$ 
\blacktriangleright F   (\text{id } \otimes \Delta ) F.
$$                 
Hence we have obtained a procedure to construct
coinduced comodules.

In Ref. \cite{ibort92} it is discussed this  construction,
called  induced  representation,
starting from a Hopf algebra $H$ and a  $H'$--comodule  for a quotient Hopf 
algebra of the first one.  The results  obtained  correspond to
unitary induced corepresentations of Hopf algebras.


\section{Spaces and representations}\label{spacesyrepresentaciones}

In this  section we present some 
situations where  geometry  and  representation theory  interact directly in order
to stress the deep relationship between them. We emphasize either those 
aspects relevant on group theory  that can be
extended to  quantum Hopf algebras or the   problems that  appear.

\subsection{Invariant integrals and unitary  representations}
\label{Invariantintegralsandunitaryrepresentations}

It is well know that the  use of sums (or integrals) and averages over a group
is very useful in group theory. Thus, now we will try to construct an invariant
integral, since  it allows  to construct a unitary module.

Let us consider the simple  case of a homogeneous $G$--space $(X,
\triangleleft, G)$, where
 $X$ as well $G$ are finite, and the  co-space 
$(F(X), \triangleright,\mathbb{C}[G])$. The dual of $F(X)$ can be
identified with the  space of formal linear combinations
$\mathbb{C}[X]$ through  the bilinear form determined by
$\langle f, x \rangle  =  f(x), \ f\in F(X), \; x \in X.$
So, a linear form on $F(X)$ will be 
$$
I ={\displaystyle \sum_{x \in X} } \alpha_x x.
$$  
If $I$ is an invariant integral then it satisfies 
$I \triangleleft h = \epsilon(h) I, \ \forall h \in  \mathbb{C}[G]$.
Taking $h=g^{-1}$ the l.h.s. of the invariance condition  is
$I \triangleleft g^{-1} = {\sum_{x \in X} } 
\alpha_x (x\triangleleft g^{-1})= { \sum_{x' \in X} } 
\alpha_{x'\triangleleft g} x'$,
and the r.h.s.
$\epsilon(g^{-1}) I = I  = { \sum_{x' \in X} } 
\alpha_{x'} x'$.
Both expressions give that $\alpha_{x'\triangleleft g} =\alpha_{x'}$. Since
$X$ is a homogeneous  space all the coefficients
 $\alpha_x$ are equal to a constant $\lambda \in \mathbb{C}$, hence
$I  =\lambda   { \sum_{x \in X} }  x$.
Let $|X|$ be the cardinal of $X$. The integral is
normalized choosing $\lambda =|X|^{-1}$, and  in this case
$I(f) =|X|^{-1}{ \sum_{x \in X} } f(x)$,
with $f$ an arbitrary function. The integral $I$ combined with the
star structure on the space of functions, i.e. $f^*(x)= \overline{f(x)}$,
allows to define an inner product $\langle f_1, f_2 \rangle =I(f_1^* f_2)$.
If the star structure on $\mathbb{C}[G]$ is given by
$g^* =g^{-1}$ we get that
$\langle g^* \triangleleft f_1, f_2 \rangle= 
\langle  f_1, g \triangleleft f_2 \rangle$.
This result can be extended by linearity to any element of
$\mathbb{C}[G]$. So, we have proved that the  module  
$(F(X),\triangleright,  \mathbb{C}[G])$ becomes unitary.

The generalization of this result is provided by the following proposition.
 \begin{propos} \label{unitarizacion}
  Let $(A, \triangleright, H)$ be a co-space, such that  
$A$ and $H$ have star  structures compatible in the sense of
$(h\triangleright a)^*= S(h)^* \triangleright a^*, \forall a\in A, \forall h\in
H$.  If $I$ is a $H$--invariant integral on $A$, then
  the star  structure defined on $H$ is compatible with the bilinear form
    $\langle a, b \rangle = I(a^*b)$.
 \end{propos}

Note that if $(A, \langle \cdot, \cdot \rangle)$
is a Hilbert space then the  representation associated to 
the $H$--module $(A, \triangleright, H)$ is unitary.

\subsection{Co-spaces and induction} \label{cospaceseinduccion}

In practice the induction algorithm has two steps: the characterization
of the carrier space $V^\uparrow =\text{Hom}_{A'}(A,V)$ as subspace  of
$\text{Hom}_{\mathbb{K}}(A,V)$ and the description of the action.
When  a Hopf  algebra is involved in this procedure  both steps can be
formulated in terms of  co-spaces.
Effectively, let us consider a subalgebra $L$ (non necessarily a  Hopf subalgebra)
of a Hopf algebra 
$H$ and a $L$--module $(V, \triangleright, L)$, and suppose that  $H$
as well as $V$ are finite dimensional for the sake of simplicity. The 
carrier space of the  induced representation 
$V^\uparrow =\text{Hom}_L(H,V)$ is the subspace of $V\otimes H^*$, whose
 elements $ \sum_i v_i \otimes \varphi_i$ verify
$  \sum_i \langle \varphi_i, lh \rangle v_i  
=\sum_i  \langle \varphi_i, h\rangle (l\triangleright v_i),
\ \forall h \in H, \; \forall l \in L$.
This condition can be expressed using the regular action of $H$ on its dual seen in
section \ref{spacesyco-spaces} (and no reference to the  elements of
$H$ is made)
$$  
\sum_i  v_i \otimes (\varphi_i \triangleleft l) 
=\sum_i (l\triangleright v_i) \otimes \varphi_i,
\qquad  \; \forall l \in L \subset H.
$$
Hence, the equivariance condition can be described in terms of the
regular co-space $(H^*,\triangleleft,H)$. The connection between the action on
the induced module and the other regular co-space is summarized 
in the following proposition.
\begin{propos}\label{proposicion2}
The  module  $(V^\uparrow, \triangleright, H)$ is a submodule of the module 
tensor product $(V, \triangleright, H) \otimes  (H^*, \triangleright, H) $, whose
second factor is the left regular $H$--co-space and the action on the first factor
is given by the counit  by means of 
$$ 
h \triangleright v  =\epsilon(h) v,
\qquad h \in H, \; v \in V.
$$
\end{propos} 
The above considerations can be formulated in a simple way when  the  module from
where one induces is unidimensional, since in this case the carrier space
$V^\uparrow$ is the subspace of $H^*$ determined by  the condition
$\varphi \triangleleft l  =\varphi, \ \forall l \in L.$
Moreover, the  induced module  is simply a submodule of the regular module 
$(H^*, \triangleright, H)$. It  suffices to consider 
Proposition \ref{proposicion2} and to take into account that  the character
$\epsilon$ of $H$ is the unit element of the tensor product. 

The final conclusion, as we mention in the introduction,   is that in the
 induction procedure  both regular co-spaces are used: $(H^*,\triangleleft,H)$
characterizes the equivariance condition, i.e., the carrier space of the induced
representation,  and by restricting  the action on $(H^*, \triangleright, H)$  the 
induced representation is obtained explicitly.

\section{Induced representations of $U_q(\overline{\fra{g}(1,1)})$}\label{q-galileo}

The quantum extended Galilei group is the Hopf algebra
$F_q(\overline{G(1,1)})$, generated by $\mu, x, t$ and $v$ with nonvanishing
commutation relations 
$$
[\mu,x]=-2a\mu, \qquad [\mu,v]= a v^2, \qquad [x,v]= 2av.
$$
The coproduct, counit and antipode are
$$
\begin{array}{c}
\begin{array}{ll}
\Delta \mu= \mu \otimes 1 + 1 \otimes \mu + v \otimes x + 
\frac{1}{2} v^2 \otimes t, &\qquad
\Delta x= x \otimes 1 + 1 \otimes x + v \otimes t, \\[2mm]
\Delta t = t \otimes 1 + 1 \otimes t ,&\qquad
\Delta v = v \otimes 1 + 1 \otimes v;
 \end{array} \\ \\[-2.5mm]
\epsilon(\mu)=\epsilon(x)=\epsilon(t)=\epsilon(v)=0 ; \\[2mm]
S(\mu)= -\mu + x - \frac{1}{2} v^2 t, \quad S(x)= -x +t v , \quad
S(t)= -t, \quad S(v)= -v.
\end{array}
$$

The dual (enveloping) algebra $U_q(\overline{\fra{g}{(1,1)}})$
\cite{Bon92} is generated by $I, P, H$ and $N$, in such a way that its
pairing with
$F_q(\overline{G(1,1)})$ is given by
$$
\langle I^p P^q H^r N^s, \mu^{p'} x^{q'} t^{r'} v^{s'}\rangle =
p! q! r! s! \; \delta^p_{p'} \delta^q_{q'} \delta^r_{r'} \delta^s_{s'}.
$$
The duality relations fix the Hopf algebra structure in
$U_q(\overline{\fra{g}{(1,1)}})$. The
nonvanishing commuting relations are
$$
[I, N]= -a e^{-2 a P} I^2, \qquad
[P, N]= - e^{-2 a P} I, \qquad
[H, N]= - \frac{1- e^{-2aP}}{2a}.
$$
 The  coproduct, counit and antipode can be written as
$$
\begin{array}{c}
\begin{array}{ll}
\Delta M= M \otimes e^{-a P} + e^{a P} \otimes M, &\qquad
\Delta P= P \otimes 1 + 1 \otimes P, \\[2mm]
\Delta H= H \otimes 1 + 1 \otimes H, &\qquad
\Delta K= K \otimes e^{-a P} + e^{a P} \otimes K;
\end{array} \\ \\[-2.5mm]
\epsilon(M)= \epsilon(P)=\epsilon(H)=\epsilon(K); \\[2mm]
S(M)= -M, \quad S(P)= -P, \quad S(H)= -H,
\quad S(K)= -K- aM;
\end{array}
$$
where  $M= e^{-a P} $ and  $K= e^{aP} N$. However, the pairing is simpler in terms
of $(I, P, H, N)$ than instead of $(M, P, H, K)$. 

After some computations (see \cite{arratia00}) one finds that :

\noindent 1).- The action of the generators of $U_q(\overline{\fra{g}(1,1)})$
on the left regular module
$(F_q(\overline{G(1,1)}),\succ,U_q(\overline{\fra{g}(1,1)}))$ is given by
$$\label{action431}
\begin{array}{cclr}
I \succ f &=& \left(1 + a \bar{v}  \frac{\partial}{\partial \mu}
e^{-2a  \frac{\partial}{\partial x}}
\right) \frac{\partial}{\partial \mu} f, &\\[2mm]
P \succ f &=& \left[ \frac{\partial}{\partial x} +
 \frac{1}{a} \ln \left(1 + a \bar{v}  \frac{\partial}{\partial \mu}
e^{-2a  \frac{\partial}{\partial x}}\right) \right] f, &\\[2mm]
H \succ f &=&  \left[ \frac{\partial}{\partial t} + \frac{1}{2a} \bar{v}
\left(1-\frac{ e^{-2a  \frac{\partial}{\partial x}}}{
1 + a \bar{v}  \frac{\partial}{\partial \mu}
e^{-2a  \frac{\partial}{\partial x}}} \right) \right] f, &\\[2mm]
N \succ f &=&  \frac{\partial}{\partial v}f &f\in
F_q(\overline{G(1,1)}).
\end{array}
$$

\noindent 2).- The action on the right regular module
$(F_q(\overline{G(1,1)}),\prec,U_q(\overline{\fra{g}(1,1)}))$  is 
$$\label{action432}
\begin{array}{ccl}
 f \prec I  &=&  \frac{\partial}{\partial \mu} f, \qquad
 f \prec P =   \frac{\partial}{\partial x} f, \qquad
 f \prec H =   \frac{\partial}{\partial t} f, \\ [0.2cm]
 f \prec N &=&  \left[ \frac{\partial}{\partial v}+
 a \bar{\mu} e^{-2a  \frac{\partial}{\partial x}}\frac{\partial^2}{\partial \mu^2} +
\bar{x} e^{-2a  \frac{\partial}{\partial x}}   \frac{\partial}{\partial \mu}+
\bar{t}  \frac{1-e^{-2a  \frac{\partial}{\partial x}} }{2a} \right] f.
\end{array}
$$

The above results allow us to construct a
family of representations of $U_q(\overline{\fra{g}(1,1)})$ coinduced by
the character $(I^p P^q H^r) \vdash 1= \alpha^p \beta^q \gamma^r$
of the Abelian subalgebra generated by $I, P$ and $H$. The carrier space
$\C^\uparrow$ is the subspace of elements of $F_q(G(1,1))$ of the form
$e^{ \alpha \mu} e^{\beta x} e^{\gamma t} \phi(v)$, whit  $\phi$  an
arbitrary function. The action can be transferred to
the space of formal power series
 $\C[[v]]$, where the action of the generators of
 $U_q(\overline{\fra{g}(1,1)})$ is
$$
\begin{array}{ccl}
 I \vdash \phi(v)& = &   \alpha (1 + a  \alpha e^{-2 a  \beta} v) \phi(v) ,
\\[2mm]
 P \vdash \phi(v)& = &  [ \beta + \frac{1}{a} \ln(1 
+ a  \alpha e^{-2 a  \beta} v)]\phi(v) ,\\[0.2cm]
H \vdash \phi(v)& = & [ \gamma+ \frac{1}{2a}(1-
 \frac{e^{-2a \beta}}{1+a \alpha e^{-2 a \beta} v})v] \phi(v) ,\\[0.2cm]
 N \vdash \phi(v) & = & \phi'(v) .
       \end{array}
$$
Induced representations of  
 $U_q(\overline{\fra{g}(1,1)})$ were previously obtained by Bone\-chi {\sl et al}
\cite{bgst98}.  We have deduced the coregular representations and constructed from
them a family of coinduced representations  including those presented in
 \cite{bgst98}. 


\section*{Acknowledgments}
\noindent 
This work has been partially supported by DGES of the Ministerio de Educaci\'on 
y Cultura de Espa\~na under Projects PB98--0719 and PB98--0360, 
and the Junta de Castilla y Le\'on (Espa\~na).
\bigskip


\end{document}